\documentclass [10pt,a4paper]{article}
\usepackage{amssymb}
\usepackage{enumerate,verbatim}
\usepackage{amsmath}
\usepackage{amsfonts,mathrsfs}
\usepackage{amsthm}
\usepackage{epsfig}
\usepackage[applemac]{inputenc}
\usepackage{pdfsync}
\usepackage{xcolor,hyperref}
\usepackage[T1]{fontenc}
\usepackage{lpic}
\usepackage[width=1.2\textwidth]{caption}

\hypersetup{backref,colorlinks=true,linkcolor=black,,citecolor=black}

\usepackage{marvosym }

\newcounter{teoremaganso}
\newcounter{coryganso}

\renewenvironment{abstract}{\small\quotation\noindent
 {\bfseries \abstractname .}}{\endquotation \par}

\newenvironment{trivlistparm}[2]{\trivlist
\item[{#1 #2}]}{\endtrivlist}

\newenvironment{prooftext}[1]{\trivlistparm{\bfseries}{#1}}{\Qed\endtrivlistparm}
\newenvironment{prova}{\trivlistparm{\bfseries}{Proof.}}{\Qed\endtrivlistparm}

\catcode`\@=11

\def\resetthefootnote{\renewcommand{\thefootnote}{\@arabic\c@footnote} }
\def\@principiremex#1{\trivlist
 \item[\hskip \labelsep{\bfseries #1\ \thetheo.}]\ignorespaces}
\def\opar@principiremex#1[#2]{\trivlist
 \item[\hskip \labelsep{\bfseries #1\ \thetheo\ (#2).}]\ignorespaces}

\newcommand{\newTHEOremrom}[2]{\newenvironment{#1}{\refstepcounter{theo}\@ifnextchar[{\opar@principiremex{#2}}
{\@principiremex{#2}}}{\qedB\endtrivlist}} \catcode`\@=12
\DeclareMathSymbol{\square}{\mathord}{AMSa}{"03}
\newcommand{\qedB}{\nopagebreak\hspace*{\fill}$\square$\par}
\newcommand{\Qed}{\nopagebreak\hspace*{\fill}{\vrule width6pt height6pt depth0pt}\par}

\newtheorem {theo} {Theorem} [section]
\newtheorem {prop} [theo] {Proposition}

\newtheorem {lem} [theo] {Lemma}
\newtheorem {bigtheo} [teoremaganso] {Theorem}
\newtheorem {bigcory} [coryganso] {Corollary}

\newTHEOremrom {defi} {Definition}
\newTHEOremrom {obs} {Remark}
\newTHEOremrom {ex} {Example}

\newcommand{\refc}[1]{\mbox{$(\ref{#1})$}}

\newcommand{\teoc}[1]{Theorem~\ref{#1}}
\newcommand{\propc}[1]{Proposition~\ref{#1}}
\newcommand{\coryc}[1]{Corollary~\ref{#1}}
\newcommand{\lemc}[1]{Lemma~\ref{#1}}
\newcommand{\defic}[1]{Definition~\ref{#1}}
\newcommand{\obsc}[1]{Remark~\ref{#1}}

\newcommand{\figc}[1]{Figure~\ref{#1}}

\newcommand{\N}{\ensuremath{\mathbb{N}}}
\newcommand{\Z}{\ensuremath{\mathbb{Z}}}
\newcommand{\R}{\ensuremath{\mathbb{R}}}

\newcommand{\T}{\ensuremath{\mathcal{T}}}

\newcommand{\C}{\ensuremath{\mathbb{C}}}

\newcommand{\D}{\ensuremath{\Theta_{\lambda}}}
\newcommand{\du}{\ensuremath{\mathcal{D}}}
\newcommand{\F}{\ensuremath{\mathcal{F}}}

\newcommand{\A}{\ensuremath{\mathcal{A}}}

\newcommand{\RP}{\ensuremath{\mathbb{RP}}}

\newcommand{\barep}{\ensuremath{\vert\hspace{-1.3mm}\varepsilon}}
\def\map#1#2#3{\mbox{${#1}\!:{#2}\longrightarrow{#3}$}}

\newcommand{\q}{\mathcal Q}
\newcommand{\V}{\mathcal V}
\newcommand{\U}{\mathcal U}
\newcommand{\UU}{\mathscr U_{1}}

\title{\bf Unfoldings of saddle-nodes and their Dulac time
\footnotetext{This work was partially supported by the FONDECYT project 1120333, 
the grants  MTM2011-26674-C02-01, MTM-2008-03437 from  MINECO/FEDER,
 UNAB10-4E-378, co-funded by ERDF ``A way to build Europe'' and by the French ANR project STAAVF.}     
}

\author{P. Marde\v si\'c, D. Mar\'{\i}n,  M. Saavedra and J. Villadelprat
\\*[.1truecm]
{\small \textsl{Institut de Math{\'e}matiques de Bourgogne, UFR Sciences et Techniques,}}
\\*[-.05truecm]
{\small \textsl{Universit\'{e} de Bourgogne, B.P. 47870, 21078 Dijon, France}}
\\*[.1truecm]
{\small \textsl{Departament de Matem\`{a}tiques, Facultat de Ci\`{e}ncies,}}
\\*[-.05truecm]
{\small \textsl{Universitat Aut\`{o}noma de Barcelona, 08193 Bellaterra, Barcelona, Spain}}
\\*[.1truecm]    
{\small \textsl{Departamento de Matem\'atica, Facultad de Ciencias F\'isicas y Matem\'aticas}}
\\*[-.05truecm]
{\small \textsl{Universidad de Concepci\'on, Barrio Universitario, Concepci\'on, Casilla 160-C, Chile}}
\\*[.1truecm]
{\small \textsl{Departament d'Enginyeria Inform\`{a}tica i Matem\`{a}tiques,
                        ETSE,}}
\\*[-.05truecm]
{\small \textsl{Universitat Rovira i Virgili, 43007
                     Tarragona, Spain}}}
                     
\date{}

\begin{document}

\maketitle

\begin{abstract} 
In this paper we study unfoldings of saddle-nodes and their Dulac time. 
By unfolding a saddle-node, saddles and nodes appear. 
In the first result (Theorem~A) we prove uniform regularity by which orbits and their derivatives arrive at a node. Uniformity is with respect to all parameters including the unfolding parameter bringing the node to a saddle-node and a parameter belonging to a space of functions. 
In the second part, we apply this first result for proving a regularity result (Theorem~B) on the Dulac time (time of Dulac map) of an unfolding of a saddle-node. 
This result is a building block in the study of bifurcations of critical periods in a neighbourhood of a polycycle. 
Finally, we apply Theorems \ref{A} and \ref{B} to the study of critical periods of the Loud family of quadratic centers and we prove that no bifurcation occurs for certain values of the parameters (Theorem C).
\end{abstract}

\section{Introduction and main results}

\begin{figure}[t]
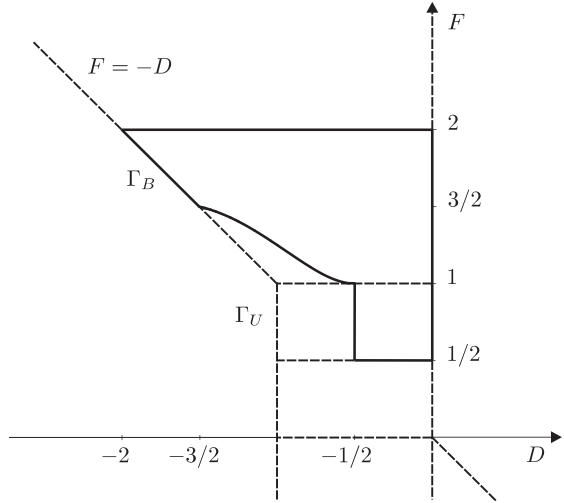

  \centering    
 \begin{lpic}[l(0mm),r(0mm),t(0mm),b(7mm)]{diagrama(0.85)}
  \end{lpic}     
\caption{Bifurcation diagram of the period function of \refc{loud} at the outer boundary according to Theorem A in~\cite{MMV2}. More precisely, $\R^2\setminus\{\Gamma_B\cup\Gamma_U\}$ corresponds to local regular values, whereas~$\Gamma_{B}$ are local bifurcation values. The results in that paper did not allow us to determine the character of the parameters in the dotted curve~$\Gamma_U$.}\label{dib3}
\end{figure}

\begin{figure}[t]
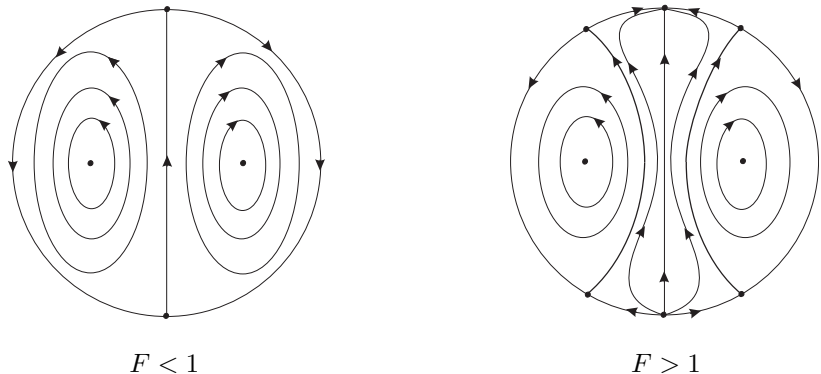

  \centering    
 \begin{lpic}[l(0mm),r(0mm),t(0mm),b(7mm)]{sella_node(2)}
  \lbl[l]{8,-3;$F<1$}
  \lbl[l]{41,-3;$F>1$}
  \end{lpic}     
\caption{Phase portrait of the Loud family \refc{loud}, with $D\in (-1,0)$ and $F>0$ in the Poincar\'eé disc, where the vertical invariant straight line is $x=1.$}\label{dib1}
\end{figure}

This paper is dedicated to the study of saddle-nodes and their unfoldings in the real plane. 
Our initial motivation comes from the study of bifurcations of critical periods of quadratic centers, but we think that our results are of more general interest.  
From the point of view of the study of the period function, the most interesting stratum of quadratic centers is given by the Loud family
\begin{equation}\label{loud}
\left\{\begin{aligned}
  & \dot u=-v+uv, \\
  & \dot v=u+Du^2+Fv^2,
\end{aligned}\right.
\end{equation}
which has a Darboux first integral. Compactifying $\R^2$ to the Poincar\'e disc, the boundary of the period annulus of the center has two connected components, the center itself and a polycycle. We call them the inner and outer boundary of the period annulus, respectively. 
In \cite{MMV2}, we described one part of the bifurcations of local critical periods from the outer boundary in this family, but many claims remained conjectural, see \figc{dib3}.
In particular, the study of the segment $\{D\in (-1,0),F=1\}$ requires a theoretical result about the local time function of such a family in a neighbourhood of a saddle-node appearing  at infinity. There, the center is bounded in the Poincar\'e disc by a symmetric polycycle and, crossing the line $F=1$, an unfolding of a saddle-node with polar factor at infinity occurs (see \figc{dib1}).
By using the  Darboux first integral and introducing an auxiliar parameter $\varepsilon=2(F-1)$,  it can be checked that this saddle-node unfolding  can be brought to the form
\begin{equation}\label{loud-equi}
\frac{1}{y U(x,y)}\bigl(x(x^{2}-\varepsilon)\partial_x-V(x)y\partial_y\bigr),
\end{equation}
by a local change of coordinates, where $y=0$ corresponds to the line at infinity in~\eqref{loud}. Taking advantage of the symmetry of the differential system~\eqref{loud}, it suffices to study half of the period and then the essential part of the period is given by the Dulac time in a neighbourhood of an unfolding of a saddle-node. In \teoc{regular}, we prove that no bifurcation of critical periods occurs for the values of the parameter in $\bigl\{D\in (-1,0)\setminus\{-\frac{1}{2}\}, F=1\bigr\}$, corresponding to saddle-nodes.

In a more general context, it is well-known that, by blowing-ups, any singularity of a vector field reduces to simple singularities or saddle-nodes.
Hence, the period function around any monodromic polycycle can always be expressed as the sum of Dulac times of saddle or saddle-node singularities, composed by their corresponding Dulac maps. Therefore, local Dulac time of saddles or saddle-nodes at finite or infinite distance can be thought of as the basic building blocks in the study of the period function near the outer boundary of a period annulus. In \cite{MMV} and \cite{MMV3}, we deal with orbitally linearizable and resonant saddles, respectively. In this paper we consider the remaining case: saddle-node singularities.

According to \cite{RT}, a family unfolding a saddle-node  is always
analytically orbitally equivalent to
\[
Y(x,y;\lambda)=P_{\lambda}(x)\partial_x+\bigl(P_{\lambda}(x)R_{0,\lambda}(x)-y(1+a(\lambda)x^{\mu})+y^2R_{2,\lambda}(x,y)\bigr)\partial_y,
\]
 where $P_{\lambda}(x)=x^{\mu+1}+\nu_{\mu-1}(\lambda)x^{\mu-1}+\cdots+\nu_1(\lambda) x+\nu_0(\lambda)$
 and $R_{i,\lambda}$ are germs of holomorphic functions.  The differential form $\omega_\lambda$ dual to the vector field $Y(x,y,\lambda)$ is
 \[
 \omega_{\lambda}=\bigl(P_{\lambda}(x)R_{0,\lambda}(x)-y(1+a(\lambda)x^{\mu})+y^2R_{2,\lambda}(x,y)\bigr)dx-P_{\lambda}(x)dy.
 \]
Functions $R_{0,\lambda}$ and $R_{2,\lambda}$ correspond to unfoldings of the Martinet-Ramis moduli of analytic classification of saddle-nodes. The vanishing of $R_{0,\lambda}$ corresponds to the existence of an analytic central manifold in the family, which by a change of coordinates can be put at $y=0$. This occurs when the saddle-node belongs to an analytic polycycle bounding the period annulus. It is also a necessary condition for Darboux local integrability, i.e. for the existence of a local meromorphic integrating factor for the differential form~$\omega_{\lambda}$. Although the vanishing of $R_{2,\lambda}$ is not a necessary condition for Darboux local integrability, it provides a sufficient condition for the Liouvillian local integrability, i.e. for the existence of a meromorphic closed differential form $\eta$ such that $d\omega=\omega\wedge\eta$. In this situation the function $\exp\!\left(\int\!\eta\right)$ is a (not necessarily meromorphic) integrating factor of $\omega$. These considerations, jointly with the example (\ref{loud-equi}) provided by the Loud family, justify to consider the model case  $R_{0,\lambda}=R_{2,\lambda}=0$ in which we have triviality of the unfolded moduli.  

One can verify that the Dulac time of an unfolding of the form
\begin{equation}\label{temp0}
\frac{1}{y^{n}U(x)}\bigl(P(x)\partial_x-V(x)y\partial_y\bigr)
\end{equation}
is an \emph{orbit} $y=y(x)$ of the saddle-node unfolding
$$P(x)\partial_x+\bigl(n V(x)y-U(x)\bigr)\partial_y$$
arriving to a nodal sector when the independent variable $x$ approaches to a specific separatrix of \eqref{temp0}. For this reason, we study in detail the regularity of the trajectories arriving to these nodal sectors, which is an interesting problem by itself. 
Conceptually, our approach to study the local Dulac time will be similar to the use of Picard-Fuchs systems in the study of Abelian integrals. 
The aim of this paper is hence two-fold.

In the first part, we deal with a class of Darboux integrable systems unfolding a saddle-node. 
We consider trajectories of the unfolding arriving to one of the \emph{nodes}. We study the regularity by which these trajectories and their derivatives arrive at the node when the unfolding parameter tends to zero (\teoc{A}). 

\begin{figure}[t]
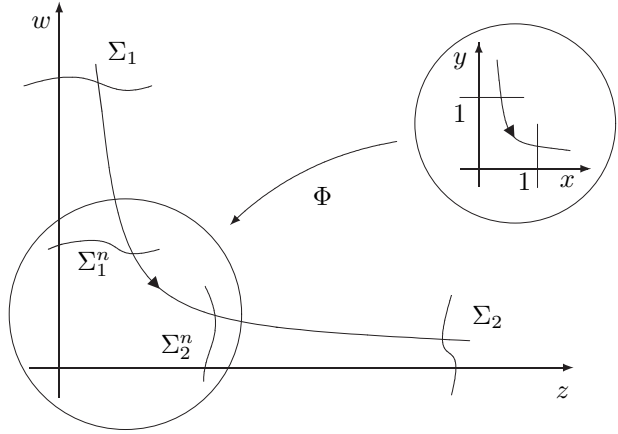

  \centering    
 \begin{lpic}[l(0mm),r(0mm),t(0mm),b(0mm)]{transversals(1)}
  \lbl[l]{3,54;$w$}
  \lbl[l]{72,5;$z$}
  \lbl[l]{72.5,33;$x$}
  \lbl[l]{67,33;$1$}
  \lbl[l]{58.5,49;$y$}
  \lbl[l]{58.5,42;$1$}
  \lbl[l]{40,31;$\Phi$}
  \lbl[l]{13,50;$\Sigma_{1}$}
  \lbl[l]{9,22;$\Sigma_{1}^{n}$}
  \lbl[l]{20,11;$\Sigma_{2}^{n}$}
  \lbl[l]{61,15;$\Sigma_{2}$}
  \end{lpic}     
\caption{Local Dulac time and normalized transverse sections.}\label{dib2}
\end{figure}

In the second part, we study the local Dulac time between normalized transverse sections (see \figc{dib2}) of the saddle part of the saddle-node unfolding. We permit a polar factor in the unfolding, which appears when singular points at infinity in the Poincar\'e disc are studied in local coordinates. We show that the local Dulac time permits an asymptotic expansion in the monomial power scale $s^r$, $r\in\Z^{+}=\N\cup\{0\}$, and the remainder term tends uniformly to zero, together with its derivatives  of any finite order.
Our approach to prove Theorem \ref{B} forces the introduction of an unbounded parameter $\lambda$ and a functional parameter~$U(x)$ in the family~\eqref{eq0} below. Finally, we particularize our main results to a simple unfolding of the saddle-node. This example contains Euler's equation, which has a divergent center manifold.
  
In order to present our main results properly, we fix $\mu\in\N$ and we consider the following unfolding of a saddle-node
\begin{equation}\label{eq0}
X=P_{\varepsilon}(x)\partial_x+\bigl(\lambda V_{a}(x)y-U(x)\bigr)\partial_y,
\end{equation}
parametrized by $(\varepsilon,a,\lambda,U)$, with $\varepsilon\approx 0,$ $a$ in an open subset $A$ of $\R^{\alpha}$, $\lambda>0$, $U\in\mathscr U$ and
\begin{itemize}
\item $P_{\varepsilon}(x)=P(x;\varepsilon)$ is an analytic function in $(x,\varepsilon)$, for $|x|\leqslant r$,  such that $P_{0}(x)$ has a zero of  order $\mu+1$ at  $x=0$;
\item $V_{a}(x)$ is analytic in $(x,a)$, for $|x|\leqslant r$, with $V_{a}(0)=1$, for all $a\in A$;
\item $\mathscr U$ is the space of series $U(x)=\sum_{j\geqslant 0}u_{j}x^{j}\in\R\{x\}$, with convergence radius greater than~$r>0$.
\end{itemize}

By rescaling, we assume that $r=1$ and $V_{a}(x)>0$, for $|x|\leqslant 1$, for all $a\in A$. We endow $\mathscr U$ with the norm $\|U\|:=\sum_{j\geqslant 0}|u_{j}|$ and with this norm it becomes a Banach space. We denote $\UU:=\{U\in\mathscr U\,:\,\|U\|\leqslant 1\}$. 

By Weierstrass preparation theorem and rescaling, we can assume that $P_{\varepsilon}(x)$ is a polynomial of  degree $\mu+1$ in $x$, with $P_{0}(x)=x^{\mu+1}$. The reason for not including the parameter~$\lambda$ into~$a$ is that it will vary in  an unbounded interval. This will play a key role in the proof of \teoc{B}. Moreover, the limit case $\lambda=\infty$ corresponds to a singular deformation or slow-fast system, which is also of independent interest. 

Notice that the singularity $(x,y)=\left(0,U(0)/\lambda\right)$ is a saddle-node of $X|_{\varepsilon=0}$, whose (real) nodal sector is contained in the half plane $x\geqslant 0$. 
We will assume that $P_{\varepsilon}(x)$ has a real root, for $\varepsilon\approx 0$.
In what follows,~$\vartheta_{\varepsilon}$ will denote the biggest real root of $P_{\varepsilon}(x)$. 
As it will be clear in a moment, our results refer to this root, and the reason for choosing this one among the others is because we can approach it from the right inside a sector that does not shrink as~$\varepsilon$ tends to zero. In the study of bifurcations, having  uniformity on the parameters is crucial and, with respect to $\varepsilon,$ this only makes sense by approaching from the right to $\vartheta_{\varepsilon}.$ Moreover, this is the only relevant situation in the study of the period function near the outer boundary of the period annulus. In the sequel, we will assume
\begin{enumerate}[{\bf (H0)}]
\item $P_{\varepsilon}'(\vartheta_{\varepsilon})>0$, for $\varepsilon\approx 0$, so that the singular point $(x,y)=\big(\vartheta_{\varepsilon},\frac{U(\vartheta_{\varepsilon})}{\lambda V_{a}(\vartheta_{\varepsilon})}\big)$ is a node of $X$.
\end{enumerate}
The polynomial $P_{\varepsilon}(x)$ need not be irreducible. We identify the two branches that contain the root $x=\vartheta_{\varepsilon}$, for $\varepsilon\geqslant 0$ and $\varepsilon\leqslant 0$, and we apply Puiseux theorem to each one, obtaining $\rho_{\pm}\in\N$
and analytic functions $\sigma_{\pm}(z)\in\R\{z\}$, such that 
\begin{equation}\label{puis}
\vartheta_{\varepsilon}=\left\{\begin{array}{rl}
\sigma_{-}\left((-\varepsilon)^{1/\rho_{-}}\right), & \text{if } \varepsilon\leqslant 0,\\[7pt]
\sigma_{+}\left((+\varepsilon)^{1/\rho_{+}}\right), & \text{if } \varepsilon\geqslant 0.
\end{array}
\right.
\end{equation}
Note that $\sigma_{\pm}(0)=0,$ because $\vartheta_{\varepsilon}$ tends to zero, as $\varepsilon\to 0.$ This gives the continuity of the function $\vartheta_{\varepsilon}.$ Note that this function in general is not analytic at $\varepsilon=0$, even though $\sigma_{-}$ and $\sigma_{+}$ are. In our first result, \teoc{A}, we treat the unfolding~\refc{eq0}, as $\varepsilon\to 0^+$, or $\varepsilon\to 0^-$. Since the substitution $\varepsilon\longmapsto-\varepsilon$ interchanges both situations, we will restrict to the case $\varepsilon\geqslant 0$ and, in what follows, when there is no risk of confusion, we will omit the subscript $+$, for the sake of shortness.

Besides the natural assumption (H0), we need to impose two technical conditions on $P_{\varepsilon}(x)=P(x;\varepsilon)$. In order to state them precisely, we introduce the function
\begin{equation}\label{defQ} 
\q(s,\varepsilon)\!:=\frac{P\bigl(s+\sigma(\varepsilon);\varepsilon^{\rho}\bigr)}{s},
\end{equation}
which is analytic at $(s,\varepsilon)=(0,0)$ and polynomial in $s$. Moreover, $\q(s,0)=s^{\mu}$ and on account of~(H0), $\q(0,\varepsilon)=\chi\,\varepsilon^{\nu}+\ldots$, with $\chi>0$, for some $\nu\in\N$. Taking this notation into account, the aforementioned assumptions are the following:
\begin{enumerate}[\bf (H1)]
\item The Newton's diagram of $\q(s,\varepsilon)$ has only one compact side (connecting the endpoints $(\mu,0)$ and $(0,\nu)$), i.e. $\q$ admits a Taylor's expansion of the form
$$\q(s,\varepsilon)=\sum\limits_{\frac{i}{\mu}+\frac{j}{\nu}\geqslant 1}q_{ij}s^{i}\varepsilon^{j}.$$ 
\item The principal $(\mu,\nu)$-quasi-homogeneous part of $\q(s,\varepsilon)$ is positive definite on the first quadrant, i.e.
$$\sum\limits_{\frac{i}{\mu}+\frac{j}{\nu}= 1}q_{ij}\sin^{i}\theta\cos^{j}\theta>0,\text{ for all }\theta\in\left[0,\frac{\pi}{2}\right].$$
\end{enumerate}
Notice that (H2) implies (H0) because $P'_{\varepsilon}(\vartheta_{\varepsilon})=\q(0,\varepsilon)$. On the other hand, (H1) implies (H2), if $\gcd(\mu,\nu)=1$.  However, the last implication does not hold in general, as the following example shows.

\begin{ex}
If $P_{\varepsilon}(x)=x((x-\varepsilon)^{2}+\varepsilon^{4})$, then $\vartheta_{\varepsilon}\equiv 0$, $\q(x,\varepsilon)=(s-\varepsilon)^{2}+\varepsilon^{4}=s^{2}-2s\,\varepsilon+\varepsilon^{2}+\varepsilon^{4}$ and $\mu=\nu=2$. One can easily show that $P_\varepsilon$ satisfies (H0) and (H1), but it does not satisfy (H2). 
\end{ex}

Let $y(x)=y(x;x_{0},y_{0},\varepsilon,a,\lambda,U)$ be the trajectory of \eqref{eq0}, i.e. the solution of the linear differential equation
\begin{equation}\label{eq_lineal}
P_{\varepsilon}(x) y'(x)=\lambda V_{a}(x)y(x)-U(x),
\end{equation}
with initial condition $y(x_{0})=y_{0}$.
We have $y(x)=D(x)\frac{y_0}{D(x_0)}+y_{L}(x),$ where
\begin{align*}  
&D(x)=D(x;\varepsilon,a,\lambda)\!:=\exp\left(\lambda\int_1^x\frac{V_a(s)}{P_{\varepsilon}(s)}ds\right)\\
\intertext{ and }
&y_{L}(x;x_{0},\varepsilon,a,\lambda,U)\!:=  D(x;\varepsilon,a,\lambda)\int_x^{x_0}\frac{U(s)}{P_{\varepsilon}(s)}\frac{ds}{D(s;\varepsilon,a,\lambda)}.
\end{align*}
Here $D(x)$ is a fundamental solution of the homogeneous equation and it coincides with the Dulac map of the saddle point $(x,y)=(\vartheta_{\varepsilon},0)$ of the vector field $P_{\varepsilon}(x)\partial_{x}-\lambda V_{a}(x)y\partial_{y}$, for $x\geqslant \vartheta_{\varepsilon}$. Moreover, $y_{L}(x)$ is the particular solution with initial condition $y_{0}=0$ and it depends linearly on $U\in\mathscr U$.
We are now in position to state our first main result where, for convenience, we use the differential operator
\begin{equation}\label{derivada}
\D =\frac{1}{\lambda}s\partial_{s}.
\end{equation}
\begin{bigtheo}\label{A}
Let us consider the saddle-node unfolding given in \refc{eq0}, with $\varepsilon\geqslant 0.$ Assume that $P_{\varepsilon}(x)$ satisfies the hypothesis~(H1) and~(H2).
Then, there exist functions $c_{j}(\varepsilon,\lambda,a,U)$, $j\in\Z^{+}$, satisfying that, for each $\ell,k\in \Z^{+}$, $\lambda_{0}>0$ and every compact set $K_{a}\subset A$, there exists $\varepsilon_{0}>0$, such that $c_{0},\ldots,c_{\ell}$ are  analytic on
$\A\!:=[0,\varepsilon_{0}]\times K_{a}\times[\lambda_{0},\infty)$ and are uniformly bounded linear operators on $\mathscr U$ and
 the following assertions hold:
\begin{enumerate}[(1)]
\item for every compact set $K_{x}\subset (0,1]$, 
the particular solution $y_{L}$ of \refc{eq_lineal} is of the form 
\[
y_{L}(s+\vartheta_{\varepsilon};x_{0},\varepsilon,a,\lambda,U)=\sum_{j=0}^{\ell}c_{j}(\varepsilon^{1/\rho},a,\lambda,U)s^{j}+s^{\ell}h_{\ell}(s;x_{0},\varepsilon,a,\lambda,U),
\]
where  $\D^{r}h_{\ell}(s)\to 0$, as $s\to 0^{+}$, uniformly on $K_{x}\!\times\A\!\times\UU$, for $r=0,1,\ldots,k$;
\item the fundamental homogeneous solution of \refc{eq_lineal} is of the form $D(s+\vartheta_{\varepsilon};\varepsilon,a,\lambda)=s^{\ell}h_{\ell}(s;\varepsilon,a,\lambda)$,
where
 $\D^{r}h_{\ell}(s)\to 0$, as $s\to 0^{+}$, uniformly on $\A$, for $r=0,1,\ldots,k$.
\end{enumerate}
\end{bigtheo}

\teoc{A} can be compared to the results \cite{R} and \cite{RT} of Rousseau and Teyssier but important differences exist. Both studies deal with unfoldings of saddle-nodes. 
Rousseau and Teyssier deal with the complex foliation, whereas our study is essentially real. They construct what they call \emph{squid sectors} on which, by a holomorphic change of coordinates, the vector field can be brought to a model and give moduli of analytic classification in terms of comparison of these normalizing coordinates and the period functions on asymptotic cycles giving the temporal part of the moduli. Their model equation is like our equation \eqref{eq0}, but with $U=0$.
The real sector $[\vartheta_{\varepsilon},x_0]$, which we study, would belong to one of their squid sectors, having the node in its boundary. Our study gives the asymptotic expansion of the solutions at the boundary of such a squid sector with a good uniform control together with all derivatives of the remainder term. It is algorithmic. We think that it cannot be obtained from the results in \cite{RT}. Note that requiring the uniform flatness property of a remainder term has proved its efficiency  in studying the cyclicity (creation of cycles) and their bifurcations from hyperbolic polycyles, see e.g.~\cite{Mou}. It is also the right condition for studying the bifurcation of critical periods from monodromic polycycles. 

A specific situation which will be useful for further applications is the following:  
\begin{equation}\label{equi}
P_{\varepsilon}(x)=x(x^{\mu}-\varepsilon)
\text{ and }U(x)=x^{m}\bar U(x),\text{ with } m\in\Z^{+}.
\end{equation}
In this case we have 
\begin{equation}\label{sigma-equi}
\vartheta_{\varepsilon}=\left\{\begin{array}{ll} 0, &\text{if }\varepsilon\leqslant 0,\\[5pt]
\varepsilon^{1/\mu},&\text{if } \varepsilon\geqslant 0.
\end{array}\right.
\end{equation}
Our next main result follows almost directly by applying twice \teoc{A}. We point out that it deals with both cases $\varepsilon\geqslant 0$ and $\varepsilon\leqslant 0$.

\begin{bigcory}\label{CorA}
Consider the saddle-node unfolding given in \eqref{eq0}, taking the functions given in~\eqref{equi} and setting~$\vartheta_{\varepsilon}$, as in~\eqref{sigma-equi}. Then there exist functions $c_{j}(\varepsilon,\lambda,a,U)$, $j\in\Z^{+}$, satisfying that, for each $\ell,k\in \Z^{+}$, $\lambda_{0}>0$ and  every compact set $K_{a}\subset A$, there exists $\varepsilon_{0}>0$ such that $c_{0},\ldots,c_{\ell}$ are  continuous on
$\A\!:=[-\varepsilon_{0},\varepsilon_{0}]\times K_{a}\times[\lambda_{0},\infty)$ and are uniformly bounded linear operators on $\mathscr U$, with $c_{j}(\varepsilon,a,\lambda,U)=0$, for $\varepsilon\leqslant 0$ and $j=0,1,\ldots,m-1;$ and such that
the following assertions hold:
\begin{enumerate}[(1)]
\item for every compact set $K_{x}\subset (0,1]$, 
the particular solution $y_{L}$ of \refc{eq_lineal} is of the form 
\[
y_{L}(s+\vartheta_{\varepsilon};x_{0},\varepsilon,a,\lambda,U)=\sum_{j=0}^{\ell}c_{j}(\varepsilon,a,\lambda,U)s^{j}+s^{\ell}h_{\ell}(s;x_{0},\varepsilon,a,\lambda,U),
\]
where $\D^{r}h_{\ell}(s)\to 0$, as $s\to 0^{+}$, uniformly on $K_{x}\!\times\A\!\times\UU$, for $r=0,1,\ldots,k$;
\item the fundamental homogeneous solution of \refc{eq_lineal} is of the form $D(s+\vartheta_{\varepsilon};\varepsilon,a,\lambda)=s^{\ell}h_{\ell}(s;\varepsilon,a,\lambda)$,
where
 $\D^{r}h_{\ell}(s)\to 0$, as $s\to 0^{+}$, uniformly on $\A$, for $r=0,1,\ldots,k$.\end{enumerate}
\end{bigcory}

It is worth to notice that the case $\mu=1$, $\varepsilon=0$, $\lambda=1$, $V_{a}\equiv 1$, $m=1$ and $\bar U\equiv 1$ in the above corollary corresponds to the classical Euler equation $x^{2}\partial_x+(y+x)\partial_y$,
having an irregular singular point at the origin and a divergent central manifold.

Now we motivate our main result concerning the second goal of the paper. The setting is the study of the period function of a family of polynomial centers in the plane. Since the Dulac time and its derivative of a singularity at finite distance tends to infinity, the interesting situation occurs when there are vertices of the polycycle bounding the period annulus that belong to the divisor at infinity obtained by desingularization. We study here the Dulac time of an unfolding of a saddle-node at infinity. Generically, the hyperbolic sector of the saddle-node belonging to the polycycle bounding the period annulus is deformed to a hyperbolic sector of a saddle point. This saddle point either remains at infinity or comes to finite distance in the unfolding. In the second situation there is a superposition of two different geometric phenomena. For this 
reason we study the first case. The simplest way to assure this situation is by requiring that the weak separatrix of the saddle-node unfolding is at infinity. Then, in suitable local coordinates, such an unfolding writes as 
\begin{equation}\label{temp}
\frac{1}{yU_{a}(x,y)}
\left(P_{\varepsilon}(x)\partial_x-V_{a}(x)y\partial_y\right).
\end{equation}
Without loss of generality, we assume that $U_{a}(x,y)\not\equiv 0$ has an absolutely convergent Taylor series at $(x,y)=(0,0)$ on $|x|,|y|\leqslant 1$, and that $V_{a}(x)$ is an analytic function on $|x|\leqslant 1$, with $V_{a}(0)>0$, for all $a\in A$. Notice that under assumption (H0), the point $(\vartheta_{\varepsilon},0)$, where~$\vartheta_{\varepsilon}$ is the biggest root of $P_{\varepsilon}(x)$, is now a \emph{saddle} of the differential system \eqref{temp}, for $\varepsilon\approx 0.$ In these local coordinates, the period annulus  is in the quadrant $y\geqslant 0$ and $x\geqslant\vartheta_{\varepsilon}$.
In the statement of our next result, $\T(s;\varepsilon,a)$ is the Dulac time of the saddle-node unfolding~\eqref{temp} between the transverse sections $\{y=1\}$ and $\{x=1\}$, i.e. it is the time that the trajectory starting at $(s+\vartheta_{\varepsilon},1)$ spends to arrive to $\{x=1\}$. We also use $\Theta=\Theta_1,$ see \refc{derivada}, for shortness.

\begin{bigtheo}\label{B}
Let us consider the Dulac time $\T(s;\varepsilon,a)$ of the saddle-node unfolding~\eqref{temp}, with $\varepsilon\geqslant 0.$ 
Assume that $P_{\varepsilon}(x)$ satisfies conditions (H1) and (H2).
Then there exist functions $c_{j}(\varepsilon,a)$, $j\in\Z^{+}$, satisfying that, for each $\ell,k\in \Z^{+}$ and every compact set $K_{a}\subset A$, there exists $\varepsilon_{0}>0$ such that $c_{0},\ldots,c_{\ell}$ are  analytic on 
$[0,\varepsilon_{0}]\times K_{a}$;
and the Dulac time can be written as $$\T(s;\varepsilon,a)=\sum_{j=0}^{\ell}c_{j}(\varepsilon^{1/\rho},a) s^{j}+s^{\ell}h_{\ell}(s;\varepsilon,a),$$
with $\Theta^{r}h_{\ell}(s)\to 0$, as $s\to 0^{+}$,
uniformly on $[0,\varepsilon_{0}]\times K_{a}$, for $r=0,1,\ldots,k$.
\end{bigtheo}

As we already did in Corollary~\ref{CorA}, we particularize the unfolding \refc{temp} considered in \teoc{B} by taking
\begin{equation}\label{equi2}
P_{\varepsilon}(x)=x(x^{\mu}-\varepsilon)
\text{ and }U_{a}(x,y)=x^{m}\bar U_{a}(x,y),
\end{equation}
where $m\in\Z^{+}$ and $a\in A$. As before, we stress that our next result deals with both cases, $\varepsilon\geqslant 0$ and $\varepsilon\leqslant 0$.

\begin{bigcory}\label{CorB}
Let us consider the Dulac time $\T(s;\varepsilon,a)$ of the saddle-node unfolding~\eqref{temp} taking the functions in \refc{equi2} and  setting~$\vartheta_{\varepsilon}$ as in~\eqref{sigma-equi}. Then there exist functions $c_{j}(\varepsilon,a)$, $j\in\Z^{+}$, satisfying that for each $\ell,k\in \Z^{+}$  and every compact set $K_{a}\subset A$, there exists $\varepsilon_{0}>0$ such that $c_{0},\ldots,c_{\ell}$ are continuous on 
$[-\varepsilon_{0},\varepsilon_{0}]\times K_{a}$;
and the Dulac time can be written as
$$
\T(s;\varepsilon,a)=\sum_{j=0}^{\ell}c_{j}(\varepsilon,a)s^{j}+s^{\ell}h_{\ell}(s;\varepsilon,a)
$$
with $\Theta^{r}h_{\ell}(s)\to 0$, as $s\to 0^{+}$,
uniformly on $[-\varepsilon_{0},\varepsilon_{0}]\times K_{a}$, for $r=0,1,\ldots,k$. Moreover, $c_{j}(\varepsilon,a)=0$, for $\varepsilon\leqslant 0$ and $j=0,1,\ldots,m-1$.
\end{bigcory}

Before introducing our main results, we mentioned that the \emph{local} Dulac time is the basic building block in the study of the period function near the polycycle at the outer boundary of the period annulus. We conclude this section by clarifying the role of \teoc{B} and \coryc{B} in this study and stating our last main result.  After blowing-up the singularities, we can decompose the period function near the polycycle as a sum of Dulac times between \emph{arbitrary} transverse sections $\Sigma_1$ and $\Sigma_2$ as it is shown in \figc{dib2}. In order to study each Dulac time, we use a diffeomorphism  that brings the unfolding of the singularity to its normal form; a saddle or a saddle-node.
In this paper, we study the saddle-node unfolding, given in~\refc{temp}. We use the normalizing diffeomorphism $\Phi$ 
to introduce two auxiliary \emph{normalized} transverse sections $\Sigma_1^n\!:=\Phi(\{y=1\})$ and $\Sigma_2^n\!:=\Phi(\{x=1\})$.  The function $\T$ in \teoc{B} and \coryc{B} is precisely this local Dulac time between $\Sigma_1^n$ and $\Sigma_2^n.$ In order to have a general result on the Dulac time between arbitrary transverse sections,
one must add to the local Dulac time the two times necessary
to go from given transverse sections to the normalized ones. For applications it is convenient to express these times in the coordinate on the source transversal and this leads to a composition
problem. The symmetry of the differential system~\refc{loud} makes this composition problem easier than the general situation and we are able to solve it with the tools developed in the present paper. More precisely, Corollaries~\ref{CorA} and~\ref{CorB}, together with a result obtained in~\cite{MMV}, enable us to answer the initial question motivating this paper. We can thus prove the following result, see also \figc{dib3}, where for a precise definition of local regular value we refer the reader to \cite{MMV2}.

\begin{bigtheo}\label{regular}
Denoting $a=(D,F)$, let $\{X_{a},a\in\R^2\}$ be the family of differential systems in \refc{loud} and consider the period function of the center at the origin. Then the parameters $a\in\bigl\{D\in (-1,0)\setminus\{-\frac{1}{2}\}, F=1\bigr\}$ are local regular values of the period function at the outer boundary. 
\end{bigtheo}

We point out that, by a result in \cite{MMV2}, the exceptional parameter $(D,F)=(-\frac{1}{2},1)$ is a local bifurcation value, as it can be seen in \figc{dib3}.

\section{Orbital results}

This section is dedicated to the proof of \teoc{A} and \coryc{A}, but first some preliminary and auxiliary results must be proved. To this end, we fix once for all $\lambda_{0}>0$ and compact subsets $K_{a}\subset A$ and 
$K_{x}\subset (0,1]$. Unless explicitly stated, we shall assume that $\varepsilon\geqslant 0$ and in the sequel we shall use
 \[
  \barep=\barep(\varepsilon)\!:=\varepsilon^{1/\rho},
 \]
where $\rho=\rho_{+}\in\N$ is the inverse of the Puiseux exponent given in \refc{puis}.  
Recall that a trajectory $y=y(x)$ of the unfolding, given in \eqref{eq0}, verifies the linear differential equation $P_{\varepsilon}(x) y'(x)=\lambda V_{a}(x)y(x)-U(x).$ Accordingly,
$$P_{\varepsilon}(s+\vartheta_\varepsilon) y'(s+\vartheta_\varepsilon)=\lambda V_{a}(s+\vartheta_\varepsilon)y(s+\vartheta_\varepsilon)-U(s+\vartheta_\varepsilon).$$
On account of $\vartheta_{\varepsilon}=\sigma(\barep)$, from the definition in \refc{defQ}, we get $\frac{P_{\varepsilon}(s+\vartheta_\varepsilon)}{s}=\q(s,\barep).$  Thus, since $\D=\frac{1}{\lambda}s\partial_s$, setting
$$\T(s,\barep)\!:=y(s+\sigma(\barep)),\;\V(s,\barep)\!:=V_{a}(s+\sigma(\barep))\text{ and }\;\U(s,\barep)\!:=\frac{1}{\lambda}U(s+\sigma(\barep)),$$
the above linear differential equation writes as
$\q\D\T=\V\T-\U.$
The idea to obtain the asymptotic expansion is to search for a formal solution $\T(s)=c_0+c_1s+c_2s^2+ \ldots $ satisfying
$$\frac{1}{\lambda}\q(s,\barep)s(c_1+2c_2s+ \ldots )=\V(s)(c_0+c_1s+c_2s^2+ \ldots )-\U(s).$$
Since $\q(s,\barep)s\big|_{s=0}=0$, evaluating in $s=0$, we get $c_0=\frac{\U(0)}{\V(0)}$.
Next step gives
\[
c_1=\lim\limits_{s\to 0}\left[\frac{1}{s}\left(\frac{\U(s)}{\V(s)}-\frac{\U(0)}{\V(0)}\right)\frac{\V(s)}{\V(s)-\frac{1}{\lambda}\q(s)}\right].
\]
We formalize this inductive procedure as follows.
\begin{defi}\label{Fc} Consider the linear
finite difference operator, acting on  functions $f(s)$ analytic at $s=0$,
given by 
\[
\bigl(\nabla f\bigr)(s)\!:=\left\{\begin{array}{cl}
\frac{f(s)-f(0)}{s}&\text{ for $s>0,$}\\[7pt]
f'(0)&\text{ for $s=0$}.
\end{array}
\right.
\]
Setting $F_{0}=\U,$ we define inductively
$F_{\ell+1}=\V_{\ell}\nabla(F_{\ell}/\V_{\ell})$,
where 
$\V_{\ell}(s)\!:=\V(s)-\frac{\ell}{\lambda}\q(s;\barep)$. Finally, define 
$$c_{\ell}\!:=\left.\frac{F_{\ell}}{\V_{\ell}}\right|_{s=0}\text{ and }\,\Sigma_{\ell}(s)\!:=\sum_{j=0}^{\ell}c_{j}s^{j}.$$
Note that, for each $\ell$, $c_{\ell}=c_{\ell}(\barep,\lambda,a,U)$ is a well defined function on $[-\varepsilon_{\ell},\varepsilon_{\ell}]\times [\lambda_{0},\infty)\times K_{a}\times\mathscr U$, for some $\varepsilon_{\ell}>0$, which may go to zero, as $\ell\to+\infty$. In the previous definitions, $\ell$ belongs to $\Z^{+}$, for convenience we define $\Sigma_{-1}:=0$.
\end{defi}

Notice that the functions $F_\ell$, $\ell\in\N$, are obtained from $F_0=\U$, by iterating a sort of finite differences operators, but conjugated by multiplication by $\V_\ell$. 

\begin{obs}\label{nablak}
Let $g(s)$ be an analytic function at $s=0.$ Then, for each $m\in\N$ and $k\in\{0,1,\ldots,m\},$ we have that $\nabla^{k}\bigl(s^{m}g(s)\bigr)=s^{m-k}g(s).$
\end{obs}

\begin{lem}\label{rec} 
$\q\D\Sigma_{\ell}=\V\Sigma_{\ell}-\U+s^{\ell+1}F_{\ell+1}$, for each $\ell\in\N\cup\{-1,0\}$.
\end{lem}

\begin{prova}
We proceed by induction on $\ell$. For $\ell=-1$, $\Sigma_{\ell}\equiv 0$ and the assertion holds. Assume now that the claim is true for $\ell-1$. Then 
\begin{align*}
\q\D\Sigma_{\ell}&=\q\D\Sigma_{\ell-1}+\q\D(c_{\ell}s^{\ell})=\V\Sigma_{\ell-1}-\U+s^{\ell}F_{\ell}+\ell\q c_{\ell}s^{\ell}=\V\Sigma_{\ell}-\U+s^{\ell}\bigl(F_{\ell}-c_{\ell}\V_{\ell}\bigr)\\&=\V\Sigma_{\ell}-\U+s^{\ell+1}F_{\ell+1},
\end{align*}
because $F_{\ell}-c_{\ell}\V_{\ell}=sF_{\ell+1}$, by definition.
\end{prova}

\begin{defi}\label{defflk} 
For each $k\in\Z^{+}$ and $d\in\{0,1\}$, we say that a real function $F(s,\barep;a,\lambda,U)$ belongs to the set $\F_{k}^{d}$, 
if it can be written as
$$F(s,\barep;a,\lambda,U)=\frac{f(s,\barep;a,\lambda,U)}{\q(s;\barep)^{k}},$$
where $\q$ verifies hypothesis (H1) and (H2) and $f$ is a function such that
\begin{enumerate}[$(a)$]
\item $f(s,\barep;a,\lambda,U)$ is analytic at $(s,\barep)=(0,0)$, for fixed $a,\lambda,U$ and it is homogeneous of degree~$d$ in $U$, more precisely, for $d=0$, it does not depend on $U$ and, for $d=1$, it depends linearly on $U$,

\item $f(s,\barep;a,\lambda,U)=\sum_{i,j\geqslant 0}f_{ij}(a,\lambda,U)s^{i}\barep^{j}$ with $f_{ij}(a,\lambda,U)\equiv 0$, for 
$\frac{i}{\mu}+\frac{j}{\nu}<k$, and
\item there exists a neighbourhood $W$ of $(s,\barep)=(0,0)$ in $\C^{2}$ such that the complex-analytic extension of $f$ in $(s,\barep)$ satisfies
\[
 \sup\left\{|f(s,\barep;a,\lambda,U)|: (s,\barep)\in W, (a,\lambda,U)\in [\lambda_{0},+\infty)\!\times K_{a}\!\times\UU\right\}<+\infty.
\]
\end{enumerate}
When we write $\frac{f}{\q^{k}}\in\F_{k}^{d}$, we will assume implicitly that $f$ satisfies conditions (a), (b) and (c).
\end{defi}

\begin{lem} \label{propflk} 
The following properties hold:
\begin{enumerate}[(1)]
\item $\F_{k}^{d}$ is stable by addition, $\F_{k}^{0}\F_{m}^{d}\subset\F_{k+m}^{d}$ and $\F_{k}^{d}\subset\F_{k+1}^{d}$;
\item $\nabla\F_{0}^{d}\subset\F_{0}^{d}$;
\item $\D\F_{k}^{d}\subset\F_{k+1}^{d}$;
\item  $\frac{s^{\mu}}{\q}\in\F_{1}^{0}$, $\V_{\ell},\frac{1}{\V_{\ell}}\in\F_{0}^{0}$ and  $\U:(s,\barep;a,\lambda,U)\longmapsto \frac{1}{\lambda}U(s+\vartheta_\varepsilon)\in\F_{0}^{1}$;
\item If $F\in\F_{k}^{d}$, then there exists a neighbourhood $W$ of $(s,\barep)=(0,0)$ in $\R^{+}\!\times\R^{+}$ such that $F$ is bounded on $W\times K_{a}\times\![\lambda_{0},\infty)\times\mathscr U_{1}$.
 \end{enumerate}
\end{lem}

\begin{prova}
Assertion (1) is straightforward. To prove assertion~(2), note first that since $\nabla$ is a linear operator it preserves the homogeneous degree $d$ on $U$. On the other hand, the condition on the Newton's diagram for belonging to $\F_0^d$ (i.e. for $k=0$) is empty. 
Let~$f$ be an element of $\F_{0}^{d}$. There exists $r_{0}>0$, such that the function $f(s,\barep;a,\lambda,U)$ is defined, for every $s\in\C$, with $|s|\leqslant r_{0}$. By applying Cauchy's integral formula to the function $s\longmapsto f(s,\barep;a,\lambda,U)$, which is analytic at $s=0$, we get
\[\nabla f(s,\barep;a,\lambda,U)=\frac{1}{2\pi i}\int_{|\zeta|=r_{0}}\frac{f(\zeta,\barep;a,\lambda,U)}{(\zeta-s)\zeta}d\zeta.\] 
If $|s|\leq r_{0}/2$, then the denominator in the integrand is bounded 
away from zero and the boundedness of  the complex analytic extension of $f$ implies the boundedness of that of $\nabla f$.\\
(3) Suppose that $\frac{f}{\q^k}\in \F_k^{d}.$ We will prove first that $\frac{\D f}{\q^{k}}\in\F_{k}^{d}$. To see this, note first that  the derivative $\D$ is a linear operator and it does not affect the condition about the Newton's diagram of $f$. On the other hand, by Cauchy's differentiation formula, we have that
$$\D f(s,\barep;a,\lambda,U)=\frac{1}{2i\pi}\int_{|\xi|=r_{0}}\frac{s\, f(\zeta,\barep;a,\lambda,U)}{\lambda(\zeta-s)^{2}}d\zeta$$
is bounded  on $W\times K_{a}\times[\lambda_{0},\infty)\times\mathscr U_{1},$ where $W$ is a neighbourhood of $(s,\barep)=(0,0)$ in $\C^2.$ In particular, taking $f=\q$ and $k=1$, we deduce that $\frac{\D\q}{\q}\in\F_{1}^{0}$. Finally, on account of $\D(\frac{f}{\q^{k}})=\frac{\D f}{\q^{k}}-k\frac{f}{\q^{k}}\frac{\D\q}{\q}$, we conclude that $\D\F_{k}^{d}\subset\F_{k+1}^{d}$, by using the assertions in~(1).\\
(4)  Obviously $\frac{s^{\mu}}{\q}\in\F^{0}_{1}$. 
Due to $V_{a}(0)=1$ and $\q(0,0)=0$, it follows that, for every $\ell\in\N$ and  $\lambda_{0}>0$, there exists a neighbourhood~$W$ of $(s,\barep)=(0,0)$ in $\C^2$ such that 
$$\frac{1}{2}\leqslant V_{a}(s+\vartheta_\varepsilon)-\frac{\ell}{\lambda}\q(s;\barep)\leqslant 2$$ 
on $W\times K_{a}\times [\lambda_{0},\infty)$. This shows that $\V_{\ell}$ and $\frac{1}{\V_{\ell}}$ belong to $\F_{0}^{0}$.
Finally, $\U\in\F^{1}_{0}$ because $\frac{U(s+\vartheta_\varepsilon)}{\lambda}$ is clearly linear in $U$ and bounded  on $W\times [\lambda_{0},\infty)\times\mathscr U_{1},$ where $W$ is a sufficiently small neighbourhood of $(s,\barep)=(0,0)$ in $\C^2$ (in this case there is no dependence on $a$).\\
(5) If $F=\frac{f}{\q^{k}}$ belongs to $\F^{d}_{k}$, then $f(r^{\nu}\sin\theta,r^{\mu}\cos\theta;\xi)=r^{k\mu\nu}\tilde f(r;\tilde\xi)\in\F_{0}^{0}$, where $\tilde\xi\!:=(\tilde a,\lambda,U)$ with $\tilde a\!:=(a,\theta)$ varying in the compact set $K_{a}\times[0,\frac{\pi}{2}]$. By \obsc{nablak}, we have that $\tilde f(r;\tilde\xi)=\nabla^{k\mu\nu}(r^{k\mu\nu}\tilde f(r;\tilde\xi))$, which belongs to $\F_{0}^{0}$ thanks to the assertion $(2)$, i.e. $\tilde f(r;\tilde\xi)$ is bounded on 
$V\!\times K_{a}\!\times\![0,\frac{\pi}{2}]\!\times\![\lambda_{0},\infty)\!\times\! \mathscr U_{1}$, where $V$ is some neighbourhood of $r=0$ in $\C.$ On the other hand, hypothesis (H1) and (H2) imply that $\q(r^{\nu}\sin\theta,r^{\mu}\cos\theta)=r^{\mu\nu}\tilde q(r,\theta)$ with $\tilde q(0,\theta)\geqslant\delta>0$, for all $\theta\in[0,\frac{\pi}{2}]$. Hence, there exists $r_{0}>0$ such that $\tilde q(r,\theta)\geqslant\delta/2$, for all $r\in[0,r_{0}]$ and $\theta\in[0,\frac{\pi}{2}]$. Accordingly, this shows that
 $F(r^{\nu}\sin\theta,r^{\mu}\cos\theta;\xi)=\frac{\tilde{f}(r,\tilde\xi)}{\tilde q(r,\theta)^k}$ is bounded, when $r\in[0,r_{0}]$, $\theta\in [0,\frac{\pi}{2}]$ and $(a,\lambda,U)\in K_{a}\times [\lambda_{0},\infty)\times \mathscr U_{1}$, as desired.
\end{prova}

The reason to require the boundedness of $f$ on $W\times K_{a}\times[\lambda_{0},\infty)\times\mathscr U_{1}$, where $W$ is a neighbourhood of the origin in $\C^2$, and not just in $\R^2,$ is illustrated by the following example.

\begin{ex}
The analytic function $f(s;\lambda)=\sin(\lambda^{2}s)$ is bounded on $\R\times[\lambda_{0},\infty)$. However  $\D f(s;\lambda)=\lambda s\cos(\lambda^{2} s)$ is not bounded on $(-s_{0},s_{0})\times[\lambda_{0},\infty)$, for any $s_{0}>0$.
Notice that, although there exists a neighbourhood of $\R\times[\lambda_{0},\infty)$ in $\C\times [\lambda_{0},\infty)$ where the analytic extension of $f$ is bounded, this function is unbounded on $U\times[\lambda_{0},\infty)$, for any neighbourhood $U$ of $s=0$ in $\C.$ Thus, $f$ does not belong to $\F_{0}^{0}$ according to Definition~\ref{defflk}. 
\end{ex}

\begin{prop}\label{*}
For each $\ell\in\Z^{+}$, $F_{\ell}\in\F_{0}^{1}$ and there exists $\varepsilon_{\ell}>0$ such that $c_{\ell}(\barep,a,\lambda,U)$ is an analytic function in $(\barep,a,\lambda)\in [0,\varepsilon_{\ell }]\times K_{a}\times[\lambda_{0},\infty)$ and a uniformly bounded linear operator on $\mathscr U$. 
\end{prop}

\begin{prova}
To prove that $F_{\ell}\in\F_{0}^{1}$, we proceed by induction on $\ell$. The case $\ell=0$ follows from assertion (4) of Lemma~\ref{propflk}, because $F_{0}=\U$. The inductive step follows easily from the recursive definition $F_{\ell+1}=\V_{\ell}\nabla(F_{\ell}\V_{\ell}^{-1})$,  by using assertions~(1), (2) and~(4) of Lemma~\ref{propflk}. Using again \lemc{propflk}, 
we deduce that $\frac{F_{\ell}}{\V_{\ell}}\in\F_{0}^{1}$, which implies that $c_{\ell}=\frac{F_{\ell}}{\V_{\ell}}\big|_{s=0}\in\F_{0}^{0}$.
By condition (c) in \defic{defflk}, we have
\[\sup\left\{|c_{\ell}(s,\barep;a,\lambda,U)|: (s,\barep)\in W, (a,\lambda,U)\in [\lambda_{0},+\infty)\!\times K_{a}\!\times\UU\right\}<+\infty,\]
for some neighborhood $W$ of $(s,\varepsilon)=(0,0)$.
 The analyticity of $c_{\ell}$ in the remaining parameters $(\barep,a,\lambda)$ follows easily from Definition~\ref{Fc} by the analyticity of 
$\frac{1}{\lambda}U(s+\sigma(\barep))$ and $V_{a}(s+\sigma(\barep))$.
\end{prova}

Next, we shall study the remainder term
\begin{equation}\label{hell}
h_{\ell}(s)\!:=\frac{\T(s)-\Sigma_{\ell}(s)}{s^{\ell}}
\end{equation}
of the asymptotic expansions in Theorem~\ref{A}.
Notice that in the case of assertion (2), $h_{\ell}(s)=s^{-\ell}\du(s)$,
where we denote 
\begin{equation}
\label{delta}
\du(s;\barep)\!:=D(s+\vartheta_{\varepsilon})=\exp\left(\lambda\int_1^{s+\vartheta_{\varepsilon}}\frac{V_a(t)}{P_{\varepsilon}(t)}dt\right).
\end{equation}

The following two lemmas give the basis of induction $k=0$ in assertion $(2)$ and $(1)$ of \teoc{A}, respectively.

\begin{lem}\label{dulac}
For each $\ell\in\Z^{+}$ and $y>0$ small enough, 
there exists $\varepsilon_{0}>0$ such  that $\frac{y^{\ell}\du(s)}{s^{\ell}\du(y)}$
and $s^{-\ell}\du(s)$ tend to zero, as $s\to 0^{+}$,
 uniformly on $[0,\varepsilon_{0}]\times K_{a}\times[\lambda_{0},\infty)$.
\end{lem}

\begin{prova} 
Note that it suffices to prove the first limit as $y<1$ is fixed and $\du(y)<1$. By definition we have that
$\frac{y^{\ell}\du(s)}{s^{\ell}\du(y)}=\exp(-B(s,y;\varepsilon,a,\lambda))$, with
\[B(s,y;\varepsilon,a,\lambda)\!:=\ell\log(s/y)+\lambda\int_{s}^{y}\frac{V_{a}(x+\vartheta_{\varepsilon})}{P_{\varepsilon}(x+\vartheta_{\varepsilon})}\,dx.\]
We must prove  that there exists $\varepsilon_{0}>0$ such that $\lim\limits_{s\to 0^{+}}B(s,y;\varepsilon,a,\lambda)=+\infty$,
uniformly on $[0,\varepsilon_{0}]\times K_{a}\times[\lambda_{0},\infty)$. By hypothesis, due to the compactness of $K_a$, there exists a positive constant $m_{1}$ such that $V_a(x)\geqslant m_{1}$, for any $x\in[s-\vartheta_{\varepsilon},y-\vartheta_{\varepsilon}]$ and $a\in K_a.$ 
Recall that $\vartheta_{\varepsilon}$ is the biggest root of $P_\varepsilon$, which tends to zero as $\varepsilon\to 0$, and that $\vartheta_{\varepsilon}=\sigma\left(\barep\right)$, with~$\sigma$ analytic at zero.
We have that $P_{\varepsilon}(x+\vartheta_\varepsilon)=x\q(x,\barep)$. Due to $\q(0,0)=0$, for every $m_{0}>0$ there exists $\varepsilon_{0}>0$ such that 
 $|\q(x,\barep)|\leqslant m_{0}$, for all $\varepsilon,x\in [0,\varepsilon_0]$. Hence,
$$B(s,y;\varepsilon,a,\lambda)=\ell\log(s/y)+\lambda\int_{s}^{y}\frac{V_{a}(x+\vartheta_{\varepsilon})}{x\q(x+\vartheta_{\varepsilon})}\,dx\geqslant\left(\lambda_{0}\frac{m_{1}}{m_{0}}-\ell\right)\log (y/s).$$
Taking $m_{0}$ small enough, we see that the right hand side tends to $+\infty$, as $s\to 0^+$.
\end{prova}

We show now the case $k=0$ in assertion (1) of \teoc{A}. To this end we write, see \eqref{hell}, 
$h_{\ell}=\frac{f_{\ell}}{g_{\ell}}$ with
\begin{equation}\label{fgell}
f_{\ell}(s)\!:=\frac{\T(s)-\Sigma_{\ell}(s)}{\du(s)}\quad \text{and}\quad g_{\ell}(s)\!:=\frac{s^{\ell}}{\du(s)},
\end{equation}
where $\du(s)$ is defined in \eqref{delta}.

\begin{lem}\label{k0}
For each $\ell\in\Z^{+}$, there exists $\varepsilon_{0}>0$ such that
$h_{\ell}(s)$ tends to zero, as $s\to 0^{+}$, uniformly on $K_{x}\times [0,\varepsilon_{0}]\times K_{a}\times[\lambda_{0},\infty)\times\mathscr U_{1}$.
\end{lem}

\begin{prova}
This will follow by applying the uniform L'H\^{o}pital's rule stated in the Appendix taking~$f_{\ell}$ and~$g_{\ell}$ as in~\eqref{fgell}. To this end, we must check that these functions verify the five conditions in \propc{ULH}. Condition (a) is obvious because ~$f_{\ell}$ and~$g_{\ell}$ are differentiable for $s>0$.
Using  that $\q\D\du=\V\du$ and applying \lemc{rec}, we deduce that
$$\D f_{\ell}=\frac{-s^{\ell+1}F_{\ell+1}}{\q\du}\quad\text{and}\quad\D g_{\ell}=\frac{-s^{\ell}\V_{\ell}}{\q\du}.$$
In particular, $\partial_{s}g_{\ell}=-\frac{\lambda s^{\ell-1}\V_{\ell}}{\q\du}<0$, for $s>0$, which shows condition $(b)$.
Moreover,  
\begin{equation*}
\frac{\partial_s f_{\ell}}{\partial_s g_{\ell}}=\frac{\D f_{\ell}}{\D g_{\ell}}=s\frac{F_{\ell+1}}{\V_{\ell}}
\end{equation*}
tends to zero, as $s\to 0^{+}$ uniformly on $[0,\varepsilon_{0}]\times K_{a}\times[\lambda_{0},\infty)\times\mathscr U_{1}$, for some $\varepsilon_0>0$ small enough.
This follows from Lemma~\ref{propflk}, taking into account that $F_{\ell+1}\in\F_0^1$, thanks to \propc{*}.
This shows that~$(c)$ and~$(d)$ are verified.
It only remains to check~$(e)$. The first part follows from \lemc{dulac}. To see the second part, we must verify that, for each fixed $s>0$ small enough, $\frac{f_{\ell}(s)}{g_{\ell}(s)}=\frac{\T-\Sigma_{\ell}}{s^{\ell}}$ is bounded on $K_{x}\times [0,\varepsilon_{0}]\times K_{a}\times[\lambda_{0},\infty)\times\mathscr U_{1}$. This follows from Proposition~\ref{*} and the expression
$$\T(s)=\du(s)\int_{s+\vartheta_\varepsilon}^{x_{0}}\frac{U(x)}{P_{\varepsilon}(x)}\frac{dx}{\du(x-\vartheta_\varepsilon)},$$
on account of $\sup\bigl\{|U(x)|;x\in[s+\vartheta_\varepsilon,x_{0}]\bigr\}\leqslant \|U\|$,
the monotonicity of the Dulac map $\du(x)$ and the inequalities $s+\vartheta_\varepsilon\leqslant x\leqslant x_{0}\leqslant 1$. We can thus apply \propc{ULH}, which shows that $h_{\ell}(s)=\frac{f_{\ell}(s)}{g_{\ell}(s)}$ tends to zero, as $s\to 0^{+}$, uniformly on $K_{x}\times [0,\varepsilon_{0}]\times K_{a}\times[\lambda_{0},\infty)\times\mathscr U_{1}$, as desired.
\end{prova}

The induction step, for assertions (1) and (2) in \teoc{A}, will be treated jointly:

\begin{prop}\label{k}
For each $\ell,k\in\Z^{+}$, there exist $v_{\ell k}\in\F_{k}^{0}$ and $w_{\ell k}\in\F_{k}^{1}$ such that 
$$\D^{k}h_{\ell}=v_{\ell k}h_{\ell+k\mu}+sw_{\ell k}.$$
\end{prop}

\begin{prova}
We proceed by induction on $k$. For $k=1$, we have that 
$$\D h_{\ell}=\frac{\D f_{\ell}}{g_{\ell}}-h_{\ell}\frac{\D g_{\ell}}{g_{\ell}}=\frac{-s F_{\ell+1}+h_{\ell}\V_{\ell}}{\q}.$$
Since $\T=\Sigma_{\ell}+s^{\ell}h_{\ell}=\Sigma_{\ell+\mu}+s^{\ell+\mu}h_{\ell+\mu}$, we get
$h_{\ell}=s^{\mu}h_{\ell+\mu}+s\Sigma_{\ell}^{\mu}$, with $\Sigma_{\ell}^{\mu}\!:=\frac{\Sigma_{\ell+\mu}-\Sigma_{\ell}}{s^{\ell+1}}=\sum_{j=0}^{\mu-1}c_{\ell+1+j}s^{j}$. Therefore,
$$\D h_{\ell}=\underbrace{\frac{s^{\mu}\V_{\ell}}{\q}}_{v_{\ell,1}}h_{\ell+\mu}+s\underbrace{\frac{\V_{\ell}\Sigma_{\ell}^{\mu}-F_{\ell+1}}{\q}}_{w_{\ell,1}}.$$
It is clear that $v_{\ell,1}\in\F_{1}^{0}$, because $\frac{ s^{\mu}}{\q}\in\F^{0}_{1}$ and $\V_{\ell}\in\F_{0}^{0}$, by Lemma~\ref{propflk}.
On the other hand, by \lemc{rec},
\begin{align*}
\V_{\ell}\Sigma_{\ell}^{\mu}&=s^{-\ell-1}\bigl(\V\Sigma_{\ell+\mu}-\V\Sigma_{\ell}\bigr)-\frac{\ell}{\lambda}\q\Sigma_{\ell}^{\mu}\\[5pt]
&=s^{-\ell-1}\big(\q\D\Sigma_{\ell+\mu}-\q\D\Sigma_{\ell}-s^{\ell+\mu+1}F_{\ell+\mu+1}+s^{\ell+1}F_{\ell+1}\bigr)-\frac{\ell}{\lambda}\q\Sigma_{\ell}^{\mu}\\[5pt]
&=\q s^{-\ell-1}\D(s^{\ell+1}\Sigma_{\ell}^{\mu})-s^{\mu}F_{\ell+\mu+1}+F_{\ell+1}-\frac{\ell}{\lambda}\q\Sigma_{\ell}^{\mu}=\q(\frac{1}{\lambda}\Sigma_{\ell}^{\mu}+\D\Sigma_{\ell}^{\mu})-s^{\mu}F_{\ell+\mu+1}+F_{\ell+1}.
\end{align*}
Hence $w_{\ell,1}=\frac{1}{\lambda}\Sigma_{\ell}^{\mu}+\D\Sigma_{\ell}^{\mu}-\frac{s^{\mu}}{\q}F_{\ell+\mu+1}\in\F_{1}^{1}$, thanks to Lemma~\ref{propflk}~and Proposition~\ref{*}.
We now complete the inductive step:
\begin{align*}
\D^{k+1}h_{\ell}&=\D(\D^{k}h_{\ell})=\D(v_{\ell k}h_{\ell+k\mu}+s w_{\ell k})=\D(v_{\ell k})h_{\ell+k\mu}+v_{\ell k}\D h_{\ell+k\mu}+s\bigl(\frac{w_{\ell k}}{\lambda}+\D w_{\ell k}\bigr)\\[5pt]
&= \D(v_{\ell k})s^{\mu}h_{\ell+(k+1)\mu}+s\D(v_{\ell k})\Sigma_{\ell+k\mu}^{\mu}+v_{\ell k}[v_{\ell+k\mu,1}h_{\ell+(k+1)\mu}+sw_{\ell+k\mu,1}]+s\bigl(\frac{w_{\ell k}}{\lambda}+\D w_{\ell k}\bigr)\\[5pt]
&=\underbrace{\bigl(\D(v_{\ell k})s^{\mu}+v_{\ell k}v_{\ell+k\mu,1}\bigr)}_{v_{\ell,k+1}}h_{\ell+(k+1)\mu}+s\underbrace{\bigl(\D(v_{\ell k})\Sigma_{\ell+k\mu}^{\mu}+v_{\ell k}w_{\ell+k\mu,1}+\frac{w_{\ell k}}{\lambda}+\D w_{\ell k}\bigr)}_{w_{\ell,k+1}}.
\end{align*}
Here $v_{\ell,k+1}\in\F_{k+1}^{0}$ and $w_{\ell,k+1}\in\F_{k+1}^{1}$ on account of the inductive hypothesis, \lemc{propflk}~and \propc{*}. This completes the proof.
\end{prova}

\begin{prooftext}{Proof of \teoc{A}.}
The coefficients $c_{j}$, for $j\in\Z^{+}$, are given in \defic{Fc}. \propc{*} shows that there exists $\varepsilon_{0}>0$, such that  $c_{0},c_{1},\ldots,c_{\ell}$ are analytic on $[0,\varepsilon_{0 }]\times\! K_{a}\!\times\![\lambda_{0},\infty)$ and are uniformly bounded linear operators on $\mathscr U$. 
By \propc{k}, we get
\[\D^{r}h_{\ell}=v_{\ell r}h_{\ell+r\mu}+s\,w_{\ell r},\]
with $v_{\ell r}$ and $w_{\ell r}$ bounded on $[0,s_{0}]\!\times\! [0,\varepsilon_{0}]\!\times\![\lambda_{0},\infty)\!\times K_{a}\!\times\!\mathscr U_{1}$, for some
$s_{0}>0$ and $\varepsilon_{0}>0$, thanks to assertion~(5) in \lemc{propflk}.
Notice that the linearity of $w_{\ell r}$ on $U$ implies that in the case (2), where $U\equiv 0$, we have $w_{\ell r}\equiv 0$. 
 We conclude that the limits, as $s$ tends to zero, in assertions (1) and (2) of \teoc{A} are zero uniformly on the corresponding parameters, using Lemma~\ref{k0} and Lemma~\ref{dulac}, respectively.
\end{prooftext}

\begin{prooftext}{Proof of Corollary~\ref{CorA}.}
It is easy to check that on the half planes $\varepsilon\geqslant 0$ and $\varepsilon\leqslant 0$ the corresponding functions $\q(s;\barep)$, given by~(\ref{sigma-equi}), satisfy hypothesis (H1) and (H2). 
To show assertion (1), we apply twice assertion (1) of \teoc{A}, with  $\rho_{-}=1$ and $\rho_{+}=\mu$, to deduce that 
\[y(s+\vartheta_{\varepsilon})=\left\{\begin{array}{ll}
\sum_{j=0}^{\ell}c_{j}^{-}\bigl(\varepsilon,a,\lambda,U\bigr)s^{j}+s^{\ell}h_{\ell}^{-}(s), &\text{ for }\varepsilon\leqslant 0,\\[7pt]
\sum_{j=0}^{\ell}c_{j}^{+}\bigl(\varepsilon^{1/\mu},a,\lambda,U\bigr)s^{j}+s^{\ell}h_{\ell}^{+}(s), &\text{ for }\varepsilon\geqslant 0,
\end{array}
\right.\]
where the functions $h_{\ell}^{\pm}(s)$, depending on the parameters $(x_{0},\varepsilon,a,\lambda,U)$, satisfy
$$\D^{r}h_{\ell}^{\pm}(s)\to 0,\quad\text{as}\quad s\to 0^{+},$$
uniformly on $K_{x}\times[-\varepsilon_{0},\varepsilon_{0}]\!\times\![\lambda_{0},\infty)\!\times\!K_{a}\times\!\mathscr U_{1}$, for $r=0,1,\ldots,k$. 
The flatness property of~$h_{\ell}^{\pm}$, together with the analyticity of $c_{j}^{\pm}\bigl(\varepsilon,a,\lambda,U\bigr)$ on $\bigl([-\varepsilon_{0},\varepsilon_{0}]\cap\{\pm\varepsilon\geqslant 0\}\bigr)\!\times\!K_{a}\!\times\![\lambda_{0},\infty)\!\times\!\mathscr U$, easily implies that,
for all $j\in\Z^{+}$, the functions
\[c_{j}(\varepsilon,a,\lambda,U)\!:=\left\{
\begin{array}{ll}
c_{j}^{-}(\varepsilon,a,\lambda,U),  & \text{for }\varepsilon\leqslant 0,\\[7pt]
c_{j}^{+}(\varepsilon^{1/\mu},a,\lambda,U),  & \text{for }\varepsilon\geqslant 0,
\end{array}
\right.\]
are continuous at $\varepsilon=0$. Moreover, the coefficients $c_{0},\ldots,c_{m-1}$ are identically zero, for $\varepsilon\leqslant 0$. 
This follows  from the fact that $U(x)=x^{m}\bar U(x)$ and $\vartheta_{\varepsilon}=0$, for $\varepsilon\leqslant 0$, by using the recursive definition of~$c_{j}$ and Remark~\ref{nablak}. Finally, the derivative properties of the function
\[h_{\ell}(s;x_{0},\varepsilon,a,\lambda,U)\!:=\left\{\begin{array}{ll}
h_{\ell}^{-}(s;x_{0},\varepsilon,a,\lambda,U),  & \text{for }\varepsilon\leqslant 0,\\[7pt]
h_{\ell}^{+}(s;x_{0},\varepsilon,a,\lambda,U),   & \text{for }\varepsilon\geqslant 0,
\end{array}
\right.\]
follow from the corresponding properties of $h_{\ell}^{\pm}$.
Assertion (2) in Corollary~\ref{CorA} is deduced from assertion (2) in \teoc{A} in a similar way.
\end{prooftext}

\section{Temporal results}

This section is dedicated to the proof of \teoc{B}, which follows by applying \teoc{A}. It will be clear now  why we need uniformity on $\lambda\in [\lambda_0,\infty)$ and $U$ varying in the Banach space~$\mathscr U$.

\begin{prooftext}{Proof of \teoc{B}.}
Consider $\ell,k\in\Z^+$ and a compact set $K_a\subset A.$ We decompose the given function $U_{a}(x,y)=\sum_{n\geqslant 1}U_{n,a}(x)y^{n-1}$, with $U_{n,a}\in\mathscr U$, for all $n\in\N$ and $a\in K_{a}$.  Since $U_{a}(x,y)$ is absolutely convergent on $|x|,|y|\leqslant 1$, the series $\sum_{n\geqslant 1}\|U_{n,a}\|y^{n}$ and all its $y\partial_{y}$ derivatives have convergence radius at least $1$. Consequently, 
\begin{equation}\label{nr}
\sum_{n\geqslant 1}n^{r}\|U_{n,a}\|<\infty,\text{ for all $r\in\Z^{+}$ and $a\in K_{a}$.}
\end{equation} 

Let $y(x;s)$ be the trajectory of the  vector field
$
 P_{\varepsilon}(x)\partial_x-V_{a}(x)y\partial_y$, 
with initial condition $y(s;s)=1.$ Note that the Dulac time in the statement is given by
\[\T(s;\barep,a)=\int_{s+\vartheta_{\varepsilon}}^{1}\frac{U_{a}\bigl(x,y(x;s)\bigr)y(x;s)}{P_{\varepsilon}(x)}\,dx=\int_{s+\vartheta_{\varepsilon}}^{1}\sum_{n\geqslant 1}\frac{U_{n,a}(x)y^{n}(x;s)}{P_{\varepsilon}(x)}\,dx.
\]
We define
\[
T_{n}(s)\!:=\int_{s}^{1}\frac{U_{n,a}(x)y^{n}(x;s)}{P_{\varepsilon}(x)}dx,
\]
whose derivative satisfies
\[
\partial_{s} T_{n}(s)=\int_{s}^{1}\frac{U_{n,a}(x)\partial_{s} y^{n}(x;s)}{P_{\varepsilon}(x)}\,dx-\frac{U_{n,a}(s)}{P_{\varepsilon}(s)}
=\frac{n V_{a}(s)}{P_{\varepsilon}(s)}T_{n}(s)-\frac{U_{n,a}(s)}{P_{\varepsilon}(s)},
\]
by using $\partial_{s} y(x;s)=y(x;s)\frac{V_{a}(s)}{P_{\varepsilon}(s)}$. This shows that $T_{n}(x)$ is the trajectory with initial condition $T_{n}(1)=0$ of the vector field obtained from~\eqref{eq0} by replacing~$U(x)$ by~$U_{n,a}(x)$, $V_{a}(x)$ by $\frac{V_{a}(x)}{V_{a}(0)}$ and $\lambda$ by $n V_{a}(0)$. We can thus apply \teoc{A}, with the given $\ell,k\in\Z^+,$ the compact set $K_a\subset A$, $\lambda_{0}\!:=\inf\{V_{a}(0):a\in K_{a}\}>0$ and $U\in\mathscr U$ to obtain the asymptotic expansion of $\T_{n}(s)\!:=T_{n}(s+\vartheta_{\varepsilon})$ at $s=0.$
So, there exists $\varepsilon_{0}>0$ and 
\[
 \T_{n}(s;\barep,a)=\sum_{j=0}^{\ell}c_{j}\bigl(\barep,a,n V_{a}(0),U_{n,a}\bigr)s^{j}+s^{\ell}h_{\ell}\bigl(s;\varepsilon,a,n V_{a}(0),U_{n,a}\bigr),
\] 
where the coefficients $c_{j}$ depend analytically on $(\barep,a,n)\in [0,\varepsilon_{0}]\!\times\! K_{a}\!\times\![1,\infty)$. Moreover,
\begin{align*}
 &\gamma_{j}\!:=\sup\left\{|c_{j}\bigl(\varepsilon,a,\lambda,U\bigr)|:(\varepsilon,a,\lambda,U)\in [0,\varepsilon_{0}]\times K_{a}\times [\lambda_{0},+\infty)\times\mathscr U_{1}\right\}<+\infty\\
 \intertext{and, for all positive $s$, small enough,}
 &M_{\ell}^{r}(s)\!:=\sup\left\{|\D^{r}h_{\ell}\bigl(s;\varepsilon,a,\lambda,U\bigr)|:(\varepsilon,a,\lambda,U)\in [0,\varepsilon_{0}]\times K_{a}\times [\lambda_{0},+\infty)\times\mathscr U_{1}\right\}<+\infty,
\end{align*}
with $M_{\ell}^{r}(s)\to 0$, as $s\to 0^{+}$, for $r=0,1,\ldots,k.$ In particular, for  $(\varepsilon,a,n)\in[0,\varepsilon_{0}]\times K_{a}\times\N$ and $r=0,1,\dots,k$ we have
\begin{equation}\label{ch} 
|c_{j}(\barep,a,n V_{a}(0),U_{n,a})|\leqslant \gamma_{j}\|U_{n,a}\|\quad\text{and}\quad |\D^{r}h_{\ell}(s;\varepsilon,a,n V_{a}(0),U_{n,a})|\leqslant M_{\ell}^{r}(s)\|U_{n,a}\|.
\end{equation}
Here, it is crucial that \teoc{A} holds for $\lambda$ unbounded and $U$ varying in the Banach space $\mathscr U$.

We define at this point the coefficients
\[
c_{j}(\barep,a)\!:=\sum_{n\geqslant 1}c_{j}(\barep,n V_{a}(0),a,U_{n,a}), \text{ for all $j\in\Z^+,$}
\]
which are well-defined because the series are uniformly convergent on $(\barep,a)\in[0,\varepsilon_{0}]\!\times\!K_{a}$ thanks 
to~(\ref{nr}), with $r=0$ and the first inequality in \refc{ch}.
In particular, these coefficients are analytic  on $(\barep,a)\in[0,\varepsilon_{0}]\!\times\!K_{a}$. On the other hand, by using the second inequality in \refc{ch}, the series
\[
h_{\ell}(s;\varepsilon,a)\!:=\sum_{n\geqslant 1}h_{\ell}(s;\varepsilon,a,n V_{a}(0),U_{n,a})
\]
is uniformly convergent on $(s,\barep,a)\in[0,s_{0}]\!\times\![0,\varepsilon_{0}]\!\times\!K_{a}$, for $s_{0}$ small enough, and  it tends to zero, as $s\to 0^{+}$, uniformly on $(\varepsilon,a)$. 
Hence, the series
\begin{align*}
\sum_{n\geqslant 1}\T_{n}(s;\barep,a)&=\sum_{n\geqslant 1}\sum_{j=0}^{\ell}c_{j}(\barep,a,n V_{a}(0),U_{n,a})s^{j}+s^{\ell}\sum_{n\geqslant 1}h_{\ell}(s;\varepsilon,a,n V_{a}(0),U_{n,a})\\&=\sum_{j=0}^{\ell}c_{j}(\barep,a)s^{j}+s^{\ell}h_{\ell}(s;\barep,a)
\end{align*}
is uniformly convergent on $(s,\barep,a)\in[0,s_{0}]\!\times\![0,\varepsilon_{0}]\!\times\!K_{a}$, because it is the sum of $\ell+2$ uniformly convergent series. For this reason, we can commute summation and integration in the following expression of the Dulac time 
\[\T(s;\barep,a)=\int_{s+\vartheta_{\varepsilon}}^{1}\sum_{n\geqslant 1}\frac{U_{n,a}(x)y^{n}(x;s)}{P_{\varepsilon}(x)}\,dx=\sum_{n\geqslant 1}\T_{n}(s;\barep,a).\]
Accordingly, $\T(s;\barep,a)=\sum_{j=0}^{\ell}c_{j}(\barep,a)s^{j}+s^{\ell}h_{\ell}(s;\barep,a)$. Finally, taking $\lambda=n V_{a}(0)$ and  \refc{derivada} into account, for $r=1,2,\ldots,k$, the series 
\[
\sum_{n\geqslant 1}\left|\Theta_{1}^{r}h_{\ell}(s;\varepsilon,a,n V_{a}(0),U_{n,a})\right|=
V_{a}^{r}(0)\sum_{n\geqslant 1}n^{r}\left|\D^{r}h_{\ell}(s;\varepsilon,a,n V_{a}(0),U_{n,a})\right|
\leqslant V_{a}(0)^{r}M_{\ell}^{r}(s)\sum_{n\geqslant 1}n^{r}\|U_{n,a}\|
\]
is uniformly convergent in $(s,\barep,a)$ and tends to zero, as $s\to 0^{+}$, uniformly on $(\barep,a)$, thanks to~\refc{nr} and~\refc{ch}.
Recall that uniform convergence of a series of functions does not imply the uniform convergence of its derivatives. However,  if $\{f_n\}$ is a sequence of functions, differentiable on $[a,b]$ and such that $\{f_n(x_0)\}$ converges for some point $x_0\in [a,b]$ and $\{f_n'\}$ converges uniformly on $[a,b]$, then $\{f_n\}$ converges uniformly on $[a,b]$ to a function~$f$ and $f'(x)=\lim_{n\to\infty}f_n'(x)$, for all $x\in [a,b]$
(see \cite[Theorem~7.17]{Rudin}). Taking this into account, we can assert that $\Theta_{1}^{r}h_{\ell}(s;\varepsilon,a)=\sum_{n\geqslant 1}\Theta_{1}^{r}h_{\ell}(s;\varepsilon,a,n V_{a}(0),U_{n,a})$ tends to zero, as $s\to 0^+$ uniformly on $[0,\varepsilon_{0}]\!\times\!K_{a}$, for all $r=0,1,\ldots,k$. This concludes the proof of the result. 
\end{prooftext}

The proof of \coryc{CorB} is completely analogous to that of \coryc{CorA}.

\section{Application to Loud's system}

\begin{prooftext}{Proof of \teoc{regular}.}
To study the passage through the unfolding of saddle-node at infinity we use the chart of~$\RP^2$ given by $(z,w)=\left(\frac{1-u}{v},\frac{1}{v}\right)$. In these coordinates the Loud differential system \refc{loud} writes as
\[
\frac{1}{w}\Bigl(z\bigl(1-F-Dz^{2}+(2D+1)zw-(D+1)w^{2}\bigr)\partial{z}+w\bigl(-F-Dz^{2}+(2D+1)zw-(D+1)w^{2}\bigr)\partial{w}\Bigr),
\]
which is a meromorphic vector field with Darboux first integral
\[
I(z,w)=\frac{w}{z}\left(1-2(F-1)\frac{g(z,w)}{z^{2}} \right) ^{\frac{1}{2(F-1)}},
\]
where $
g(z,w)\!:= {\frac { \left( 2\,D+1\right) }{  \left( 2\,F-1 \right)D }}zw-{\frac 
{ \left( D+1 \right) }{  2  F D}}{w}^{2}-\frac{1}{2D}.$ One can verify that the local change of coordinates given by
\begin{equation}\label{canvi}
 \left\{x=\frac{z}{\sqrt{g(z,w)}},y=\frac{w}{\sqrt{g(z,w)}}\right\}
\end{equation}
brings the above vector field to~\refc{temp}
\[
\frac{1}{yU_{a}(x,y)}\bigl((x^{2}-\varepsilon)x\partial{x}-(2F-x^{2})y\partial{y}\bigr),
\]
with $a:=(D,F)$, 
$U_{a}(x,y)\!:=\left({\frac { \left( 2D+1 \right)}{2(2F-1)}} xy-{\frac {
\left( D+1 \right) }{4F}}{y}^{2}-\frac{D}{2}\right)^{\frac{-1}{2}}$
and particularizing
$\varepsilon\!:=2(F-1)$. 
In these local coordinates, the period annulus  is in the quadrant $y\geqslant 0$ and $x\geqslant\vartheta_{\varepsilon}$, where  $\vartheta_{\varepsilon}$ is given by \refc{sigma-equi} with $\mu=2$.
Working on a compact subset $K_{a}$ of $\bigl\{D\in (-1,0),F>\frac{1}{2}\}$, we see that $U_{a}(x,y)$ has an absolutely convergent Taylor series at $(x,y)=(0,0)$ on $|x|,|y|\leqslant r$ for some $r>0$ depending only on $K_{a}$. By rescaling the local coordinates, we can assume that $r=1$. Let  $\Phi$ be the local diffeomorphism such that $(z,w)=\Phi(x,y),$ i.e. the one obtained by inverting~\refc{canvi}.

Since the Loud system \refc{loud} is invariant by the symmetry $(u,v)\longmapsto (u,-v)$, half of the period function is the Dulac time  $T$   of the singular point at infinity between  transverse sections $\Sigma_{1}:=\{v=0,u\approx 1\}$ and $\Sigma_{2}:=\{v=0,u \approx -\infty\}$. We decompose it in three parts. Let $T_{2}(s)$ be the local Dulac time between the normalized transverse sections $\Sigma_{1}^{n}\!:=\Phi(\{y=1\})$ and $\Sigma_{2}^{n}\!:=\Phi(\{x=1\})$ starting at the point $\Phi(s+\vartheta_{\varepsilon},0)$, see Figure~\ref{dib2}.
Let $T_{1}(s)$ be the time that the trajectory starting at $\Sigma_{1}$  spends to arrive to the point $\Phi(s+\vartheta_{\varepsilon},1)$ in $\Sigma_{1}^{n}$ and let $T_{3}(s)$ be the time that the trajectory  starting at the point $\Phi(1,s)$ spends to arrive to $\Sigma_{2}$. Finally, let $\du(s)$ be the  Dulac map between  $\Sigma_{1}^{n}$ and~$\Sigma_{2}^{n}$, i.e. $\du(s)$ is defined so that the trajectory starting at $\Phi(s+\vartheta_{\varepsilon},1)$ intersects $\Sigma_{2}^{n}$ at $\Phi(1,\du(s))$. By construction, see \figc{dib2}, we have $T(s)=T_{1}(s)+T_{2}(s)+T_{3}(\du(s))$. We now examine the asymptotic expansion of each piece. To this end, we denote by $\mathcal I(A)$ the space of functions $h(s;a)$, analytic on $s\in (0,s_{0}),$ such that $h(s;a)$ and $s\partial_{s}h(s;a)$ tend to zero, as $s\to 0^{+}$ uniformly, for $a$ varying in any compact subset of $A.$ Observe that this space is stable with respect to addition and multiplication. We apply Corollary~\ref{B} (with $\ell=k=1$) to obtain 
$\varepsilon_0>0$ and a uniform asymptotic expansion
$T_{2}(s)=c_{2,0}+c_{2,1}s+sh_{2}(s)$, with $h_{2}\in\mathcal I(A)$ and $c_{2,j}$ continuous on $A\!:=K_{a}\cap\{F\in (1-\varepsilon_{0},1+\varepsilon_{0})\}.$ As we already remarked just before \teoc{A}, the graph of the Dulac map $\du(s)$ is a trajectory of the nodal vector field
$(x^{2}-\varepsilon)x\partial{x}+2F(1-\frac{x^{2}}{2F})y\partial{y}$.
Hence, by applying assertion (2) of \coryc{A} with $\lambda=2F\geqslant 1$, we deduce that $\du(s)=s h_{0}(s)$, with $h_{0}\in\mathcal I(A).$
On the other hand, the time function $T_{3}(s)$ is analytic in $s$, whereas $T_{1}(s)$ is an analytic function on $s$ composed with the continuous function $(s,\varepsilon)\mapsto s+\vartheta_{\varepsilon}$. Accordingly,
they can be written as $T_{i}(s)=c_{i,0}+c_{i,1}s+sh_{i}(s)$, with $h_{i}\in\mathcal I(A)$ and $c_{i,0}$ and $c_{i,1}$ continuous on $A$, for $i=1,3.$ Note that $T_{3}\bigl(\du(s)\bigr)=c_{3,0}+s\hat h_{3}(s)$, with $\hat h_{3}(s)\!:=c_{3,1}h_{0}(s)+h_{0}(s)h_{3}\bigl(sh_{0}(s)\bigr)$ and it can be easily checked that $\hat h_{3}\in\mathcal I(A).$ Summing up the three terms we obtain that the period function of the Loud system is of the form 
\[
 P(s;D,F)=2T(s;D,F)=c_{0}(D,F)+c_{1}(D,F)s+sh(s;D,F),
\]
with $h\!:=2(h_{1}+h_{2}+\hat h_{3})\in\mathcal I(A)$ and $c_{i}$ continuous on $A.$
On the other hand, restricting to $A\cap\{F\in (\frac{1}{2},1)\}$ the singularity at 
$(x,y)=(\vartheta_{\varepsilon},0)$ is a linearizable saddle and we can apply \cite[Proposition 5.2]{MMV} to obtain the asymptotic expansion of the period function working with a different parametrization, say $\hat s$. The two parametrizations differ by composition with a diffeomorphism $\hat s=r(s)$ such that $r(0)=0$ and $r'(0)=\alpha(D,F)\neq 0$, for $F=1.$ In this other parametrization the coefficient $\hat c_{1}(D,F)$ of $\hat s$ is explicitly calculated
\[
 \hat c_{1}(D,F)=\frac{\sqrt{\pi}(2D+1)}{\sqrt{F(D+1)^{3}}}\frac{\Gamma((3F-1)/(2F))}{\Gamma((4F-1)/(2F))}.
\]
Since $c_{1}(D,F)=\alpha(D,F)\hat c_{1}(D,F)$ and one can verify that
\[
 \lim_{F\to 1^{-}}\hat c_{1}(D,F)=\frac{2(2D+1)}{(D+1)^{\frac{3}{2}}},
 \]
it follows that $c_{1}(D,1)\neq 0$, for $D\in (-1,0)\setminus\{-\frac{1}{2}\}.$ On account of the continuity of $c_{1}$ and $h\in\mathcal I(A)$, we conclude that 
\[P'(s;D,F)=c_{1}(D,F)+h(s;D,F)+sh'(s;D,F)\neq 0,\] in a neighbourhood of any point $(s,D,F)=(0,D_{0},1)$ in $(0,1)\times (-1,0)\times (\frac{1}{2},2)$, with $D_{0}\in(-1,0)\setminus\{-\frac{1}{2}\}.$
This concludes the proof of the result.
\end{prooftext}

\section*{Appendix}

In our approach to the proof of \teoc{A}, the use of L'H\^{o}pital's rule with uniformity in the parameters is fundamental. We have not found such a version in the literature. For this reason we present here the precise statement that we need  together with a proof of it.

\begin{prop}\label{ULH}
Consider two functions $\map{f_{\nu},g_{\nu}}{(a,b)}{\R}$ depending on a parameter $\nu$ belonging to an arbitrary topological space $\Lambda$ and verifying the following:
\begin{enumerate}[$(a)$]
\item  $f_{\nu}$ and $g_{\nu}$ are differentiable on $(a,b)$,
\item $g_{\nu}'(x)\neq 0$, for all $x\in(a,b)$ and $\nu\in\Lambda$,
\item for all $\nu\in\Lambda,$ there exists $L_{\nu}\in\R$ such that 
         $\lim\limits_{x\to a^{+}}\frac{f_{\nu}'(x)}{g_{\nu}'(x)}=L_{\nu}$ uniformly on $\nu\in\Lambda$,
\item $\sup\left\{|L_{\nu}|;\nu\in\Lambda\right\}<+\infty$,
\item there exists $c\in(a,b)$  such that, for each $y\in(a,c),$ we have that 
        $\lim\limits_{x\to a^{+}}\left|\frac{g_{\nu}(x)}{g_{\nu}(y)}\right|=+\infty$, uniformly on 
        $\nu\in\Lambda$ and $\sup\left\{\left|\frac{f_{\nu}(y)}{g_{\nu}(y)}\right|;\nu\in\Lambda\right\}<+\infty$.
\end{enumerate}
Then $\lim\limits_{x\to a^{+}}\frac{f_{\nu}(x)}{g_{\nu}(x)}=L_{\nu}$, uniformly on $\nu\in\Lambda$.
\end{prop}

\begin{prova}
For a given $\varepsilon>0$, we must find $\delta>0$ such that, if $x\in (a,a+\delta)$, then $\left|\frac{f_{\nu}(x)}{g_{\nu}(x)}-L_{\nu}\right|<\varepsilon$, for all $\nu\in\Lambda.$ Let us take $\varepsilon_1\!:=\min\left(\frac{\varepsilon}{M+3},1\right)$, where $M\!:=\sup\limits_{\nu\in\Lambda}|L_{\nu}|,$ which is well defined thanks to assumption~$(d)$. From~$(c)$ there exists $\delta_1>0$ such that, if $c\in (a,a+\delta_1),$ then $\left|\frac{f'_{\nu}(c)}{g'_{\nu}(c)}-L_{\nu}\right|<\varepsilon_1$, for all $\nu\in\Lambda$. Let us fix any $y\in (a,a+\delta_1)$. By the Mean Value Theorem, for each $x\in (a,y)$, there exists $c=c_{x,y,\nu}\in (x,y)\subset (a,a+\delta_1)$ such that $\frac{f_{\nu}(x)-f_{\nu}(y)}{g_{\nu}(x)-g_{\nu}(y)}=\frac{f_{\nu}'(c)}{g_{\nu}'(c)}$. Accordingly,
\begin{equation}\label{ulh_eq1}
\left|
\frac{\frac{f_{\nu}(x)}{g_{\nu}(x)}-\frac{f_{\nu}(y)}{g_{\nu}(x)}}{1-\frac{g_{\nu}(y)}{g_{\nu}(x)}}-L_{\nu}
\right|=
\left|
\frac{f_{\nu}'(c)}{g_{\nu}'(c)}-L_{\nu}
\right|<\varepsilon_1.
\end{equation}
On the other hand, the assumption $(e)$ guarantees that there exists $z_{y}\in (a,y)$ such that, if $x\in (a,z_{y}),$ then\begin{equation}\label{ulh_eq2}
\left| \frac{f_{\nu}(y)}{g_{\nu}(x)}\right|<\varepsilon_1\,
\mbox{ and }
\left|\frac{g_{\nu}(y)}{g_{\nu}(x)}\right|<\varepsilon_1\,,
\mbox{ for all $\nu\in\Lambda.$}
\end{equation}
Here, we also used that $\frac{f_{\nu}(y)}{g_{\nu}(x)}=\frac{f_{\nu}(y)}{g_{\nu}(y)}\frac{g_{\nu}(y)}{g_{\nu}(x)}$ tends to zero uniformly on $\nu\in\Lambda$, as $x\to a^+$. 
Note then that $\left|(L_{\nu}\pm\varepsilon_1)\frac{g_{\nu}(y)}{g_{\nu}(x)}\right|
< \bigl(|L_{\nu}|+\varepsilon_1\bigr)\varepsilon_1$, and thus
\begin{equation}\label{ulh_eq3}
-\bigl(|L_{\nu}|+\varepsilon_1\bigr)\varepsilon_1<(L_{\nu}\pm\varepsilon_1)\frac{g_{\nu}(y)}{g_{\nu}(x)}<
\bigl(|L_{\nu}|+\varepsilon_1\bigr)\varepsilon_1.
\end{equation}
The second inequality in \refc{ulh_eq2} shows in particular that $1-\frac{g_{\nu}(y)}{g_{\nu}(x)}>0$, because $\varepsilon_1<1,$ so that, from~\refc{ulh_eq1}, we get
\[
  (-\varepsilon_1+L_{\nu})\left(1-\frac{g_{\nu}(y)}{g_{\nu}(x)}\right)+\frac{f_{\nu}(y)}{g_{\nu}(x)}
  <\frac{f_{\nu}(x)}{g_{\nu}(x)}<
  (\varepsilon_1+L_{\nu})\left(1-\frac{g_{\nu}(y)}{g_{\nu}(x)}\right)+\frac{f_{\nu}(y)}{g_{\nu}(x)}.
\]
Therefore,
\[
 -\varepsilon_1-(L_{\nu}-\varepsilon_1)\frac{g_{\nu}(y)}{g_{\nu}(x)}+\frac{f_{\nu}(y)}{g_{\nu}(x)}
  <\frac{f_{\nu}(x)}{g_{\nu}(x)}-L_{\nu}<
 \varepsilon_1-(L_{\nu}+\varepsilon_1)\frac{g_{\nu}(y)}{g_{\nu}(x)}+\frac{f_{\nu}(y)}{g_{\nu}(x)} .
\]
From this, on account of \refc{ulh_eq3}  and the first inequality in \refc{ulh_eq2}, we get that
 \[
 -2\varepsilon_1-(|L_{\nu}|+\varepsilon_1)\varepsilon_1<\frac{f_{\nu}(x)}{g_{\nu}(x)}-L_{\nu}<
 2\varepsilon_1+(|L_{\nu}|+\varepsilon_1)\varepsilon_1.
 \]
Accordingly,
\[
 \left|\frac{f_{\nu}(x)}{g_{\nu}(x)}-L_{\nu}\right|< \varepsilon_1(2+|L_{\nu}|+\varepsilon_1)<\varepsilon_1(3+|L_{\nu}|)
 <\varepsilon_1(3+M)<\varepsilon,
\]
as desired, and so, taking $\delta=z_{y}-a$, the result follows. 
\end{prova}

\bibliographystyle{plain}

\end{document}